\documentclass[12pt, leqno]{article}
\usepackage{amssymb}
\input epsf
\def\q \m#1#2{{\raise1pt\hbox{$#1$}\kern-1pt\big/
               \kern-1pt\raise-1pt\hbox{$#2$}}}

\def\bN{{\rm \bf N}}
\def\bR{{\rm \bf R}}
\def\bZ{{\rm \bf Z}}
\def\bQ{{\rm \bf Q}}
\def\bC{{\rm {\bf C}}}
\def\ch{{\rm  ch}}
\def\bH{{\rm \bf H}}
\def\sR{{ \rm \scriptsize  \bf R}}

\def\sC{{ \rm \scriptsize  \bf C}}

\newfam\msbfam
\font\twelmsb=msbm10 at 12pt
\font\tenmsb=msbm10 at 10 pt
\font\sevenmsb=msbm10 at 7pt
\font\fivemsb=msbm10 at 5pt

\textfont\msbfam=\twelmsb
\textfont\msbfam=\tenmsb
\scriptfont\msbfam=\sevenmsb
\scriptscriptfont\msbfam=\fivemsb

\newtheorem{thm}{Theorem}[section]
\newtheorem{lemma}{Lemma}[section]
\newtheorem{defn}{Definition}[section]
 
 \newtheorem{prop}{Proposition}[section]
 \newtheorem{cor}{Corollary}[section]

\newcommand{\C}{{\cal C}}
\newcommand{\tr}{{\rm Tr}}

\begin{document}

\renewcommand{\theequation}{\thesection.\arabic{equation}}
\setcounter{equation}{0}

\centerline{\Large {\bf On Family Rigidity Theorems I}}

\vskip 10mm  
\centerline{\bf Kefeng   LIU and Xiaonan MA}
\vskip 8mm

{\bf Abstract.} In this paper, we first prove a local family version of the 
Atiyah-Bott-Segal-Singer Lefschetz fixed point formula, then we extend 
the famous Witten's rigidity Theorems to the family case. Several family vanishing theorems for elliptic genera are also proved.\\

{\bf 0 Introduction}. Let $M,B$ be two compact smooth manifolds, 
and $\pi: M\to B$ be a  submersion with compact fibre $X$.
Assume that a compact Lie group $G$ acts fiberwise on $M$, that is the action
preserves each fiber of $\pi$. Let $P$ be a family of elliptic operators
along the fiber $X$, commuting with the action $G$. Then the family index 
of $P$ is 
\begin{eqnarray} 
{\rm Ind} (P) = {\rm Ker } P - {\rm Coker } P \in K_{G} (B).
\end{eqnarray}

Note that ${\rm Ind} (P)$ is a virtual $G$-representation. 
Let $\ch_g ({\rm Ind} (P))$ with $g\in G$ be the equivariant 
Chern character of ${\rm Ind} (P)$ evaluated at $g$. 

In this paper we will first prove a family fixed point formula which expresses 
$\ch_g ({\rm Ind} (P))$ in terms of the geometric data on the fixed points
 $X^g$ of the fiber of $\pi$. The by applying this formula, we 
generalize the Witten rigidity theorems and several vanishing theorems proved
 in [{\bf Liu4}] for 
elliptic genera to the family case. 

A family elliptic operator $P$ is called 
{\em rigid on the equivariant Chern character level }
with respect to this $S^1$-action, if $\ch_g ({\rm Ind} (P)) \in H^* (B)$
 is independent of $g \in S^1$. When the base $B$ is a point, we recover the 
classical rigidity and vanishing theorems. When $B$ is a manifold, we get 
many nontrivial higher order rigidity and vanishing theorems by taking the 
coefficients of certain expansion of $\ch_g$. For the history of the Witten 
rigidity theorems, we refer the reader to [{\bf BT}], [{\bf K}], [{\bf L}], 
[{\bf H}], [{\bf Liu1}] and [{\bf Liu2}]. The family vanishing theorems 
generalize those vanishing theorems in  [{\bf Liu4}], which in turn give 
us many higher order vanishing theorems in the family case. 
In a forthcoming paper, we will extend our results to general loop group 
representations and prove much more general family vanishing theorems which generalize the results in [{\bf Liu4}]. We believe there should be some applications of our results to topology and geometry, which we hope to report on a later occasion.

This paper is organized as follows. In Section 1, we prove the  equivariant 
family index theorem. In Section 2, we prove the family rigidity 
theorem. In the last part  of Section 2, 
 motivated by the family rigidity theorem, we state a conjecture.
In Section 3, we generalize the family rigidity theorem to the nonzero 
anomaly case. As corollaries, we derive several vanishing theorems.

\newpage

\section{ \normalsize Equivariant family index theorem}
\setcounter{equation}{0}

The purpose of this section is to prove an equivariant family index theorem. 
As pointed out by Atiyah and Singer, we can introduce equivariant 
families by proceeding as in [{\bf AS1, 2}]. Here we will  prove it directly by using the local index theory as developed by Bismut.

This section is organized as follows: 
In Section 1.1, we state our main result, Theorem 1.1.
 In Section 1.2, by using the local index theory, we prove Theorem 1.1.

\subsection{\normalsize The index bundle}

Let $M,B$ be two compact manifolds, and $\pi: M\to B$ be a fibration with 
compact fibre $X$, and assume that $\dim X = 2 k$. 
Let $TX$ denote the relative tangent bundle.
Let $W$ be a  complex vector bundle on $M$ and $h^W$ be an Hermitian metric on $W$. 

Let $h^{TX}$ be a Riemannian metric on $TX$ and $\nabla^{TX}$ be the 
corresponding Levi-Civita  connection on $TX$ along the fibre $X$. 
Then the Clifford bundle 
$C(TX)$ is the bundle of Clifford algebras over $M$ whose fibre at 
$x\in M$ is the Clifford algebra $C(T_x X)$ of $(TX, h^{TX})$. 

We assume that the bundle $TX$ is spin on $M$.
 Let $\Delta= \Delta^+ \oplus \Delta^-$ be the spinor bundle of $TX$. 
We denote by $c(\cdot)$ the Clifford action of $C(TX)$ on $\Delta$. 

Let $\nabla$ be the connection on $ \Delta$ induced by $\nabla^{TX}$. 
Let $\nabla^W$ be  a Hermitian connection on $(W, h^W)$ with curvature $R^W$. 
Let $\nabla^{\Delta \otimes W}$ be the connection on 
$\Delta {\otimes} W$ along the fibre $X$:
\begin{eqnarray}
\nabla^{\Delta \otimes W}= \nabla \otimes 1 + 1 \otimes \nabla^W.
\end{eqnarray}

For $b\in B$, we denote by $E_b, E_{\pm,b}$  the set of $\C^{\infty}$-sections 
of  $\Delta\otimes W$, $\Delta_{\pm}\otimes W$  over the fiber $X_b$.
We regard the $E_b$ as the fibre of a smooth ${\bf Z}_2$-graded infinite
 dimensional vector bundle over $B$. Smooth sections of $E$ over $B$ 
will be identified to smooth sections of $\Delta\otimes W$ over $M$. 

Let $\{ e_i\}$ be an orthonormal basis of $(TX, h^{TX})$,
let $\{ e^i\} $ be its dual basis.
\begin{defn} Define the twisted Dirac operator to be
\begin{eqnarray}\begin{array}{l}
D^X = \sum_ic(e_i) \nabla^{\Delta\otimes W}_{e_i}.
\end{array}\end{eqnarray}
\end{defn}
Then $D^X$ is a family Dirac operator which acts fiberwise on the 
fibers of $\pi$. For $b\in B$, $D^X_b$ denote the restriction of $D^X$ to the fibre $E_b$. $D^X$ interchanges $E_+$ and $E_-$. Let $D^X_{\pm}$ be the  
restrictions of $D^X$ to $E_{\pm}$.
By Atiyah and Singer [{\bf AS2}], the difference bundle over $B$ 
\begin{eqnarray}
{\rm Ind} (D) = {\rm Ker} D_{+,b}-  {\rm Ker} D_{-,b}.
\end{eqnarray}
is well-defined in the $K$-group $K(B)$.

Now, let $G$ be a compact Lie group which acts fiberwise  on $M$.
We will consider that $G$ acts as identity on $B$. Without loss of generality we can assume that $G$ acts on $(TX, h^{TX})$ isometrically. We also assume that the action of $G$ lifts to $\Delta$ and $W$, and that the $G$-action commutes with $\nabla^W$.

In this case, we know that ${\rm Ind} (D^X)\in K_G(B)$. Now we start to give a proof of a local family fixed point formula which extends [{\bf AS2}, Proposition 2.2]

\begin{prop} There exist $V_j \in \widehat{G}$ with $j=1, \cdots, r$, 
a finite number of sections $(s_{i_j+ 1}, \cdots$, $s_{i_{j+1}})$ with 
$ i_{j+1}- i_{j}= \dim V_j$ of $\C ^{\infty} (B, E_-)$ such that 
we can find  a basis $\{ e_{j,l}\} $ of $V_j$, under which the map $\overline{D}_{+,b} : \C^{\infty} (B, E_{+,b} )\oplus 
\oplus_{j=1}^r V_j \to \C^{\infty} (B, E_{-,b} )$ given by 
\begin{eqnarray}
\overline{D}^X_{+,b} (s + \Sigma_{j,l} \lambda_{j,l}e_{j,l}) =D_+^X s 
+ \Sigma_{j,l} \lambda_{j,l} s_{i_j + l}
\end{eqnarray}
is $G$-equivariant and surjective. The vector spaces 
${\rm Ker }\overline{D}^X_{+,b}$ form a $G$-vector bundle 
${\rm Ker } \overline{D}_+^X$
 on $B$, and the element 
$[{\rm Ker }\overline{D}^X_+]-\oplus_{j=1}^r V_j \in K_G(B)$
 depends only on $D^X$ and not on the choice  of $\{ V_j\}$
 and the sections $\{ s_i \}$.
\end{prop}

$Proof$: Given $b_0\in B$, we can find $a>0$ and a ball $U(b_0)\subset B$
 around $b_0$, such that for any $b\in U(b_0)$, $a$ is not an eigenvalue of 
$D^{X,2}_b$. 

Let $E^{[0, a[}_b=E^{[0, a[}_{+b} \oplus E^{[0, a[}_{-b} $ be the direct sum of the eigenspaces of $D^{X,2}_b$
associated to the eigenvalues $\lambda \in [0,a[$. 
By [{\bf{ BeGeV}}, Proposition 9.10], $E^{[0, a[}$ forms a finite dimensional 
sub-bundle $E^{[0, a[} \subset E$ over $U(b_0)$. 
Clearly, $E^{[0, a[}$ is a $G$-vector bundle on $U(b_0)$. 
By [{\bf S}, Proposition 2.2 ],
we have  an isomorphism  of vector bundles on $B$
\begin{eqnarray}
E^{[0, a[}= \oplus_{V\in \widehat{G}} {\rm Hom}_G (V, E^{[0, a[}) \otimes V,
\end{eqnarray}where $\widehat{G}$ denotes the space of all irreducible representations of $G$.
We can also find $t_{i,k}\in \C^{\infty}(U(b_0), {\rm Hom}_G (V, E^{[0, a[}_-))$
such that for $b \in U(b_0)$, the elements $t_{i,l}$ form a basis 
of ${\rm Hom}_G (V, E^{[0, a[}_-)_b$.  Let $\{ e_{i,l}\}$ be a basis of $V_i$. Then we can choose the sections 
$t_{i,k}e_{i,l}\in  \C ^{\infty} (B, E^{[0, a[}_-)$ to be our $s_i$'s. This proves the first part of the proposition locally.

The global version now follows easily by extending the above local sections of 
$\C^{\infty} (U(b_0),$ $ E_-)$ together with a use of the partition of unity argument. This is essentially the same as the proof of 
[{\bf AS2}, Proposition 2.2]. \hfill $\blacksquare$\\

By [{\bf S}, Proposition 2.2], we have 
\begin{eqnarray}
{\rm Ind} (D^X) = \oplus_{V\in \widehat{G}} {\rm Hom}_G (V, {\rm Ind} (D^X))
\otimes V
\end{eqnarray}
and ${\rm Hom}_G (V, {\rm Ind} (D^X)) \in K(B)$.
We denote by $({\rm Ind} (D^X))^G\in K(B)$ the G-invariant part of 
${\rm Ind} (D^X)$.

By composing the action of $G$ and the Chern character of 
${\rm Hom}_G (V,{\rm Ind} (D^X))$, we get 
the equivariant Chern character 
$\ch_g ({\rm Ind} (D^X)) \in H^* (B)$. 

\begin{defn} We say that the operator $D^X$  is rigid on the equivariant Chern 
character level, if $\ch_g ({\rm Ind} (D^X))$ is constant on $g\in G$.
More generally, we say $D^X$ is rigid on the equivariant $K$-Theory level, if 
${\rm Ind} (D^X) = ({\rm Ind} (D^X))^G$.  
\end{defn}
In the rest of this paper, when we say $D^X$ is rigid, we always mean 
$D^X$ is rigid on the equivariant Chern character level.\\

Now let us calculate the equivariant Chern character 
$\ch_g ({\rm Ind} (D^X))$ in terms of the fixed point data of $g$.

Let $T^H M$ be a $G$-equivariant sub-bundle of $TM$ such that 
\begin{eqnarray}
TM = T^H M \oplus TX.
\end{eqnarray}
Let $P^{TX}$ denote the projection from $TM$ to $TX$.
If $U\in T B$, let $U^H$ denote the lift of $U$
 in $T^H M$, so that 
$\pi_* U^H = U$.

Let $h^{TB}$ be a Riemannian metric on $B$, and assume that $W$ has the 
Riemannian metric 
$h^{TM}= h^{TX} \oplus \pi^* h^{TB}$. Note that our final results will be 
independent of $g^{TB}$. Let $\nabla^{TM}$, $\nabla^{TB}$ 
 denote the corresponding Levi-Civita connections  on $M$ and $B$. 
Put $\nabla^{TX}= P^{TX} \nabla^{TM}$ which is a connection on $TX$. 
As shown in [{\bf B1}, Theorem 1.9], 
$\nabla^{TX}$ is independent of the choice of $h^{TB}$. 
Now the connection $\nabla^{TX}$ is well defined on $TX$ and on $M$.
Let $R^{TX}$ be the corresponding curvature.
We denote by $\nabla$ and $ \nabla^{\Delta \otimes W}$  the corresponding 
connections on 
$\Delta$ and $ \Delta {\otimes} W$ induced by $\nabla^{TX}$ and $\nabla^W$. 

Take $g \in G$ and set
\begin{eqnarray}
M^g = \{ x\in M, g x = x\}.
\end{eqnarray}
Then  $\pi: M^g \to B$ is a fibration with compact fibre $X^g$.
By [{\bf BeGeV}, Proposition 6.14], $TX^g$ is  naturally oriented in $M^g$.

Let $N$ denote the normal bundle of $M^g$, then $N= TX/TX^g$. We denote the
differential of $g$ by $dg$ which gives a bundle isometry $dg: N\to N$.
  Since $g$ lies in a compact abelian 
Lie group, we know that there is an orthogonal decomposition $N= N(\pi) 
\bigoplus
\bigoplus_{0<\theta < \pi} N(\theta)$, where $d g_{|N(\pi)} = - {\rm id}$, 
and for each $\theta, 0<\theta < \pi$, $N(\theta)$ is a complex vector
 bundle on which $d g$ acts by multiplication by $e^{i \theta}$,
and $\dim N(\pi)$ is even. So $ N(\pi)$  also is naturally oriented.

 As the Levi-Civita connection $\nabla^{TM}$ preserves the decomposition 
$TM= TM^g \oplus_{0< \theta \leq \pi} N(\theta)$,
the connection $\nabla^{TX}$ also preserves the decomposition 
$TX= TX^g \oplus _{0< \theta \leq \pi} N(\theta)$ on $M^g$.
Let $\nabla^{TX^g}$, $\nabla^N$, $ \nabla^{N(\theta)}$
 be the corresponding induced connections on $TX^g$, $N$ and $N(\theta)$, and 
et $R^{TX^g}$, $R^N$, $ R^{N(\theta)}$ be the corresponding curvatures.
 Here we consider $N(\theta)$ as a real vector bundle.
Then we have the decompositions:
\begin{eqnarray}\begin{array}{l}
R^{TX}= R^{TX^g}\oplus R^N, \qquad R^N= \oplus_{\theta} R^{N(\theta)}.
\end{array}\end{eqnarray}

\begin{defn} For $ 0<\theta \leq \pi$, we write
\begin{eqnarray}\begin{array}{l}
\ch_g(W, \nabla^W)= \tr \Big [g \exp({ -R^{W}\over 2\pi i})\Big ],\\

\displaystyle{
\widehat{A}(TX^g, \nabla^{TX^g}) 
= {\det}^{1/2} \Big ({{i \over 4 \pi} R^{TX^g} 
\over \sinh({i \over 4 \pi} R^{TX^g})}\Big )},\\
\displaystyle{
\widehat{A}_{\theta}(N(\theta), \nabla^{N(\theta)})= 
{ 1 \over     i^{{1 \over 2} \dim N(\theta)} {\det}^{1 / 2} \Big ( 
1 - g \exp({i \over 2 \pi} R^{N(\theta)}) \Big )}    }
\end{array}\end{eqnarray}
\end{defn}
Let $\ch_g(W)$, $\widehat{A}(TX^g)$, 
 $\widehat{A}_{\theta}(N(\theta))$ denote the corresponding 
cohomology classes  on $M^g$.
 
If we denote by $\{x_j, -x_j \}$ $(j=1, \cdots, l)$  the Chern roots of 
$N(\theta)$, $TX^g$ such that $\Pi x_j$ define the orientation of 
$N(\theta)$ and $TX^g$, then 
\begin{eqnarray}\begin{array}{l}
\displaystyle{
\widehat{A}(TX^g)=\Pi_j { x_j \over 2}/ \sinh ({ x_j \over 2}),}\\
\displaystyle{
\widehat{A}_{\theta} (N(\theta))= 2^{-l} \Pi_{j=1}^l
{1 \over \sinh{1 \over 2} (x_j + i \theta)} 
= \Pi_{j=1}^l { e^{{1 \over 2} (x_j + i \theta)} 
\over e^{x_j + i \theta} -1}  .  }
\end{array}\end{eqnarray}

We denote by $\pi_* : H^*(M^g) \to H^*(B)$ the intergration along the fibre $X^g$.
\begin{thm} We have the following identity in $H^*(B)$:
\begin{eqnarray}
\ch_g({\rm Ind} (D^X))= \pi_* \left \{ \Pi_{0< \theta \leq \pi}
\widehat{A}_{\theta} (N(\theta))\widehat{A} (TX^g) \ch_g (W) \right \}.
\end{eqnarray}
\end{thm}

\subsection{\normalsize A heat kernel proof of Theorem 1.1}

As Atiyah and Singer indicated in the end of [{\bf AS2}], we can proceed 
as in [{\bf AS1,2}] to introduce an equivariant family,
 and then to find a formula for the equivariant Chern character of the index 
bundle. Here, we will use a different approach by combining the local 
relative index theory and the equivariant technique to give a direct 
proof of the local version of Theorem 1.1.

We denote by  ${^0 \nabla} = \nabla^{TX}\oplus  \pi^* \nabla^{TB}$ 
the  connection on $TM$. Let $S=\nabla^{TM}- {^0 \nabla}$.
By [{\bf B1}, Theorem 1.9 ], $\left \langle S(.).,.\right \rangle_{h^{TM}}$ 
is a tensor independent of $h^{TB}$. For $U\in T^H M$, 
we define a horizontal $1$-form $k$ on $M$ by
\begin{eqnarray}
k (U)= \sum_i \left \langle S(U)e_i , e_i\right \rangle.
\end{eqnarray}

\begin{defn} Let $\nabla^E$ denote the connection on $E$ such that if 
$U\in T B$ and $s$ is a smooth section of $E$ over $B$, then
\begin{eqnarray}
 \nabla _U^E s = \nabla _{U^H}^{ \Delta \otimes  W} s.
\end{eqnarray}
\end{defn}
If $U,V$ are smooth vector fields on $B$, we write
\begin{eqnarray}
T(U^H, V^H)= - P^{TX}[U^H, V^H]
\end{eqnarray}
which is a tensor.

Let $f_1,\cdots, f_m$ be a basis of $T B$,
 and $f^1, \cdots, f^m$ be the dual basis. Define
\begin{eqnarray}
c(T)= \frac{1}{2} \Sigma_{\alpha, \beta}
 f^{\alpha}f^{\beta}c(T(f^H_{\alpha},f^H_{\beta})).
\end{eqnarray} 

\begin{defn} For $t>0$, let $A_t$ be the Bismut superconnection constructed 
in [{\bf B1}, \S 3],
\begin{eqnarray}
A_t = \sqrt{t} D^X + (\nabla^E + {1 \over 2} k) - {1 \over 4 \sqrt{t} } c(T).
\end{eqnarray}
\end{defn}
It is clear that $A_t$ is also $G$-invariant.

Let $dv_X$ denote the Riemannian volume element on the fiber $X$.
Let $\Phi  $ be the scaling homomorphism from $\Lambda (T^* B)$ into itself :
$ \omega \to  (2 \pi i) ^{- (\deg \omega )/2} \omega$.

\begin{thm} For any $t>0$, the form $\Phi  \tr_s [g \exp(-A_t^2)]$ is closed 
and that its cohomology class is independent of $t$ and represents $\ch_g({\rm Ind} (D^X))$ in cohomology.
\end{thm}

$Proof$: Just proceed as in [{\bf B1}, \S 2(d)].
\hfill $\blacksquare$

\begin{thm} We have the following identity
\begin{eqnarray}\begin{array}{l}
\lim_{t\to 0}\Phi \tr_s[g \exp(-A_t^2)]=\\
\hspace {10mm} 
\int_{X^g} \widehat{A}(TX^g, \nabla^{TX^g}) 
\Pi_{0< \theta \leq \pi}\widehat{A}_{\theta}(N(\theta), \nabla^{N(\theta)})
\ch_g(W, \nabla^W).
\end{array}\end{eqnarray}
\end{thm}

$Proof$: If $A$ is a smooth section of $T^* X \otimes \Lambda (T^* B) \otimes {\rm End} (\Delta \otimes W)$, we use the notation
\begin{eqnarray}
(\nabla^{\Delta  \otimes W}_{e_i} + A(e_i))^2 = \sum_{i=1}^{2k} (\nabla^{\Delta  \otimes W}_{e_i} + A(e_i))^2 - \nabla^{\Delta  \otimes W}_{\Sigma_{i=1}^{2k}
\nabla ^{TX}_{ e_i} e_i}- A(\Sigma_{i=1}^{2k} \nabla ^{TX}_{ e_i}e_i).
  \nonumber
 \end{eqnarray}
Let $\nabla'_t$ be the  connection on $\Lambda (T^* B)
 \widehat{\otimes} \Delta
{\otimes} W$ on the fibre $X$.
\begin{eqnarray}
\nabla'_t = \nabla^{\Delta  \otimes W} 
+ \frac{1}{2\sqrt{t}} \left \langle S(.)e_j, f_{\alpha}^H \right \rangle 
 c(e_j)f^{\alpha}
+ \frac{1}{4t}\left \langle S(.) f_{\alpha}^H, f_{\beta}^H \right \rangle  
f^{\alpha} f^{\beta}
\end{eqnarray}
Let $K^X$ denote the scalar  curvature of the fiber $(X, h^{TX})$.
By the Lichnerowicz formula [{\bf B1}, Theorem 3.5], we get
\begin{eqnarray}\begin{array}{l}
\displaystyle{
A^2 _t= - t (\nabla'_{t,e_i})^2 
+ \frac{t}{4} K^X +  \frac{t}{2}c(e_i)c(e_j) R^W(e_i,e_j)  }\\
\displaystyle{\hspace*{10mm}
+ \sqrt{t}c(e_i)  f^{\alpha} R^W(e_i,f_{\alpha}^H)
+ \frac{1 }{2}f^{\alpha} f^{\beta} R^W (f_{\alpha}^H,  f_{\beta}^H).}
\end{array}\end{eqnarray}
Let $P_u(x,x',b) (b\in B,x,x'\in X_b)$ be the smooth kernel 
associated to $\exp(-A^2_t)$
with respect to $dv_X (x')$. Then
\begin{eqnarray}
\Phi \tr_s[g\exp(-A^2_t)]= \int_{X}\Phi \tr_s[gP_t(g^{-1}x,x,b)]dv_X (x).
\end{eqnarray}
 By using standard estimates on the heat kernel, for $b\in B$, we can reduce the 
problem of calculating the limit  of (1.21) when 
$t  \rightarrow 0$ to an open neighbourhood 
${\cal U}_{\varepsilon}$ of $X^g_b$ in $X_b$. Using normal geodesic coordinates
to $X^g_b$ in $X_b$, we will identify ${\cal U}_{\varepsilon}$ to an 
$\varepsilon$-neighbourhood of $X^g$ in $N_{X^g/X}$. We know that, if 
$(x,z)\in  N_{X^g/X}$ with $x\in X^g$, then  
\begin{eqnarray}
g^{-1} (x,z)= (x,g^{-1} z).
\end{eqnarray}

Let $dv_{X^g}(x), dv_{N_{X^g/X,x}}$ with $x\in X^g$ be the corresponding 
volume forms on $TX^g$ and $N_{X^g/X}$ induced by $h^{TX}$.
Let $k(x,z)(x\in X^g, z\in N_{X^g/X}, |z|< \varepsilon)$ be defined by
\begin{eqnarray}
dv_X = k(x,z) dv_{X^g}(x) dv_{N_{X^g/X}}(z).
\end{eqnarray}
Then it is clear that 
\begin{eqnarray}
k(x,0)=1. \nonumber
\end{eqnarray}

By the  discussion following (1.21), (1.23), we get
\begin{eqnarray}\begin{array}{l}
\lim_{t \to 0}\Phi \tr_s[g\exp(-A^2_t)] = \lim_{t \to 0}
\int_{{\cal U}_{\frac{\varepsilon}{8}}}\Phi \tr_s \Big [gP_t(g^{-1}x,x) \Big ]
dv_X (x)\\
= \lim_{t \to 0}\int_{x\in X^g} \int_{\stackrel{|Y|\leq \varepsilon /8,} 
{Y\in N_{X^g/X}}}\Phi
 \tr_s \Big [gP_t \Big  (g^{-1}(x,Y),(x,Y)\Big )\Big ] k(x,Y)\\
\hspace*{35mm} \quad  dv_{X^g} (x) 
dv_{N_{X^g/X}}(Y).
\end{array}\end{eqnarray}

By taking $x_0 \in X^g_b$ and using the finite propagation speed as in 
[{\bf B2}, $\S$ 11b)], 
we may assume that in  $X_b$ we have the identification 
$(TX)_{x_0}\simeq \bR^{2k}$  with $0\in \bR^{2k}$
representing $x_0$ and that the extended fibration over $\bR^{2k}$ 
coincides with the given fibration restricted to $B(0,\varepsilon)$.

Take any vector $Y\in \mbox{\bf R}^{2k}$. We can trivialize 
 $   \Lambda  (T^* B) \widehat{\otimes}  \Delta
{\otimes} W$  by  parallel transport along the curve 
$u \rightarrow uY$ with respect to  $\nabla'_t$.

Let $\rho(Y)$ be a $\cal C^{\infty}$-function over $\mbox{\bf R}^{2k}$ 
which is equal to 1 if
$|Y| \leq \frac{\varepsilon}{4}$, and  equal to 
0 if $|Y| \geq \frac{\varepsilon}{2}$. 
Let $\Delta^{TX}$ be the ordinary Laplacian operator on $(TX)_{x_0}$.
Let $H_{x_0}$ be the vector space of smooth sections of the bundle 
$(\Lambda  (T^* B) \widehat{\otimes}  \Delta 
{\otimes} W)_{x_0}$ over $(T X)_{x_0}$.
 For $t>0$, let $L^1_t$ be the operator  acting on $H_{x_0}$: 
\begin{eqnarray}
L^1_t = (1- \rho^2 (Y))(- t\Delta^{TX}) +  \rho^2 (Y) A^2_t.
\end{eqnarray}
For $t>0$, $s\in  H_{x_0}$, we write 
\begin{eqnarray}\begin{array}{l}
F_{t} s(Y)= s( \frac{Y}{\sqrt{t}}),\\
L^2 _t = F^{-1}_{t}L^1_t F_{t}.
\end{array}\end{eqnarray}
Let $\{ e_1,\cdots, e_{2l'}\}$ be an orthonormal  basis of 
$(T X^g)_{x_0}$, and  let
$\{ e_{2l'+1}, \cdots, e_{2k}\}$ be an orthonormal basis of $N_{X^g/X, x_0}$. 
Let $L^3_t$ be the operator obtained from $L^2_t$ by replacing the Clifford 
variables $c(e_j)$ with $1 \leq j \leq 2l'$ by the operators 
$\frac{e^j}{\sqrt{t}}- \sqrt{t} i_{e_j}$.

Let $P^i_t(Y,Y')$ with $ Y,Y'\in (T X)_{x_0}$ and $|Y'|< \frac{\varepsilon}{4}$, $i=1,2,3$ be the smooth kernel associated to $ \exp(-L^i_t)$ with respect to 
 the volume element $dv_{TX_{x_0}}(Y')$. 
By using the finite propagation speed method, there exist $c,C>0$, such that 
for $Y\in N_{X^g/X, x_0}$,
 $|Y|\leq \frac{\varepsilon}{8}$ and $t\in ]0,1]$, we have 
\begin{eqnarray}
 \Big | P_t(g^{-1}Y,Y)k(x_0, Y)  -P^1_t(g^{-1}Y,Y) \Big |
\leq c \exp(-{C \over t^2}).
\end{eqnarray}

For $\alpha \in \mbox{\bf C} (e^j, i_{e_j})_{(1\leq j \leq 2l')}$,  let 
$[\alpha]^{\mbox{\scriptsize  max}} \in  \mbox{\bf C} $ be the coefficient
of $e^1 \wedge \cdots \wedge e^{2l'}$ in the expansion of $\alpha$.
Then as in [{\bf B2}, Proposition 11.12], if $Y\in N_{X^g/X}$
\begin{eqnarray}
\qquad \mbox{Tr}_s\Big [g P^1_t(g^{-1}Y,Y)\Big ]= (-2i)^{{1 \over 2} \dim X^g}
t^{- {1 \over 2}\dim N_{X^g/X}}
\mbox{Tr}_s \Big [g P^3_t
\Big (\frac{g^{-1}Y}{\sqrt{t}},\frac{Y}{\sqrt{t}}\Big ) \Big ]
^{\mbox{\scriptsize  max}}.
\end{eqnarray}

Let $R^{TX}_{|M^g},  R^W_{|M^g}, \cdots $ be the corresponding restrictions 
of $R^{TX}, R^W,\cdots $ to $M^g$.
Let $\nabla_{e_j}$ be the ordinary differentiation operator on $(T X)_{x_0}$
in the direction $e_j$. By [{\bf ABoP}, Proposition 3.7], and  (1.20), we have,
 as $t \rightarrow 0$,
\begin{eqnarray}
L^3_t \rightarrow L^3_0 = = -  \sum_{j=1}^{2k} \Big (\nabla_{e_j} 
+ \frac{1}{4}\left \langle R^{TX}_{|M^g}Y,e_j\right \rangle  
 \Big )^2 + R^W_{|M^g}.
\end{eqnarray}

By proceeding as in [{\bf B2}, \S 11g)- \S 11i)], we obtain the following: 
there exist some constants $\gamma >0, c>0,C>0, r\in  \mbox{\bf N} $ such that for
 $t\in ]0,1]$ and  $Y,Y'\in (T X)_{x_0}$, we have 
\begin{eqnarray}\begin{array}{l}
 \quad \Big |P^3_t(Y,Y') \Big |\leq c (1+|Y|+|Y'|)^r \exp(-C|Y-Y'|^2), \\
 \Big |(P^3_t- P^3_0) (Y,Y') \Big | \leq c t^{\gamma}  (1+|Y|+|Y'|)^r
 \exp(-C|Y-Y'|^2).
 \end{array}\end{eqnarray}
From  (1.28) and (1.30), we get 
\begin{eqnarray}\begin{array}{l}
\lim_{t \rightarrow 0} \int_{\stackrel{|Y|\leq \frac{\varepsilon}{8}}
{Y\in N_{X^g/X}}}
\Phi \mbox{Tr}_s[g P^1_t(g^{-1}Y,Y)] dv_{N_{X^g/X}}(Y)\\
= \lim_{t \rightarrow 0} \int_{\stackrel{|Y|\leq 
\varepsilon /8 \sqrt{t}} {Y\in N_{X^g/X}}}
(-2i)^{ {1 \over 2} \dim X^g} \Phi
 \tr_s  [gP^3_t   (g^{-1} Y,Y  )]
dv_{N_{X^g/X}}(Y) \\
= \int_{N_{X^g/X}}  (-2i)^{ {1 \over 2} \dim X^g} 
\Phi \mbox{Tr}_s[g P^3_0(g^{-1}Y,Y)]^{\mbox{\scriptsize  max}} 
dv_{N_{X^g/X}}(Y).
 \end{array}\end{eqnarray}

Now we define 
\begin{eqnarray}
A= - \sum_{j=1}^{2k}  \Big (\nabla_{e_j} 
+ \frac{1}{4}\left \langle R^{TX}_{|M^g}Y,e_j\right \rangle  
 \Big )^2 
\end{eqnarray}
By Mehler's formula [{\bf G}], the smooth kernel $q(Y,Y')$ for $Y,Y'\in TX$, 
associated to $\exp(-A)$ is given by
\begin{eqnarray}\begin{array}{l}
\displaystyle{
q(Y, Y')= (4\pi)^{-k} {\det}^{1/2}\Big ({ R^{TX} /2 
\over \sinh( R^{TX}/2)} \Big )
\exp \Big \{ - {1\over 4}\left \langle 
{ R^{TX}/2 \over \tanh ( R^{TX}/2)}  Y, Y\right \rangle }\\
\displaystyle{\hspace*{10mm} 
 - {1\over 4}\left \langle
{ R^{TX}/2 \over \tanh ( R^{TX}/2)}  Y', Y'\right \rangle 
+ {1\over 2} \left \langle 
{ R^{TX}/2 \over \sinh ( R^{TX}/2)} e^{R^{TX}/2} Y, Y'\right \rangle \Big \}}
\end{array}\end{eqnarray}
From (1.9) and (1.33), we deduce, for $Y\in N_{X^g/X}$,
\begin{eqnarray}\begin{array}{l}
\displaystyle{q(g^{-1}Y, Y) = (4\pi)^{-k} {\det}^{1/2} \Big ({ R^{TX} /2
\over \sinh( R^{TX}/2)} \Big ) }\\

\displaystyle{\hspace*{20mm} 
\exp \Big \{ - {1 \over 2} \left \langle { R^{N}/2 \over \sinh ( R^{N}/2)} 
 \Big (\cosh (R^{N}/2)- e^{R^N/2} g^{-1} \Big ) Y,Y\right \rangle  \Big \}.}
\end{array}\end{eqnarray}
On the other hand, for $Y \in N(\theta)$, we have 
 \begin{eqnarray}
\qquad \left \langle { R^{N} e^{R^N/2}\over 2 \sinh ( R^{N}/2)}  g^{-1}Y,
Y\right \rangle  
= \left \langle { R^{N}/2 \over \sinh ( R^{N}/2)} 
{1 \over 2} ( e^{R^N/2} g^{-1} +e^{-R^N/2} g) Y,
Y\right \rangle .
\end{eqnarray}
It is easy to see that 
\begin{eqnarray}
\qquad \cosh(R^{N}/2)-{1 \over 2} ( e^{R^N/2} g^{-1} +e^{-R^N/2} g)
=  {1 \over 2} (1 - g^{-1}) ( e^{R^N/2}-   e^{-R^N/2}g).
\end{eqnarray}

From (1.9), (1.34)- (1.36),  we get
\begin{eqnarray}\begin{array}{l}
\int_{N_{X^g/X}}  q(g^{-1} Y, Y) d v_{N_{X^g/X}}(Y)
 =  (4 \pi)^{- {1\over 2} \dim X^g}\\
\displaystyle{\hspace*{10mm} 
{\det}^{1/2} \Big ({  R^{TX^g}/2 \over \sinh (R^{TX^g}/2 )} \Big )
 \Big [ {\det}^{1/2} (1- g_{|N}^{-1}) 
{\det}^{1/2} \Big (1 - g e^{-R^N}\Big )\Big ]^{-1}.}
\end{array}\end{eqnarray}

We may and will assume that on the basis $\{e_{m}\}_{2l'+1 \leq m \leq 2k}$,
 the matrix of $g$ has diagonal blocks
$$\left [ \begin{array}{l}
\cos(\theta_j) \quad -\sin(\theta_j) \\
\sin(\theta_j) \quad \cos(\theta_j)
\end{array} \right ], 0< \theta_j \leq \pi.$$
Then one verifies easily that the action of $g$ on $\Delta$ is given by
\begin{eqnarray}
g= \Pi_{l'+1\leq j\leq k} \Big (\cos(\theta_j /2)
+ \sin(\theta_j/2) c(e_{2j-1}) c(e_{2j}) \Big ).
\end{eqnarray}
By (1.29) and  (1.38),  we know that
\begin{eqnarray}
\qquad \quad  \mbox{Tr}_s \Big [g P^3_0(g^{-1}Y,Y) \Big ]
=  \Pi_{l'+1\leq j\leq k} \Big (-2i \sin(\theta_j/2)\Big )
 \tr  \Big [g \exp (-R^W_{|M^g}) \Big ] q(g^{-1} Y,Y).
\end{eqnarray}

From (1.24), (1.27), (1.31), (1.37) and  (1.39), 
we finally arrive at the wanted formula (1.18). 
\hfill $\blacksquare$\\

By Theorems 1.2 and  1.3, we now have the complete proof of Theorem 1.1.
\hfill $\blacksquare$\\

\section{ \normalsize Family Rigidity Theorem}
\setcounter{equation}{0}

This section is organized as follows. In Section 2.1,
 we state our main theorem of the paper: the family rigidity theorem.
 In Section 2.2, we prove it 
by using the  equivariant family index theorem and the modular invariance.
In Section 2.3,  motivated by the family Witten rigidity theorem, 
we state a conjecture about a $K$-theory level rigidity theorem for elliptic genera.

Throughout this section, we use the  notations of Section 1, and take $G= S^1$.

\subsection{ \normalsize Family rigidity theorem}

Let $\pi: M\to B$ be a  fibration of compact manifolds with fiber $X$ and $\dim X= 2k$. We assume that the $S^1$ acts fiberwise on $M$,
and $TX$ has an $S^1$-equivariant spin structure. As in [{\bf AH}], 
by lifting to the double cover of $S^1$, the second condition is always 
satisfied. Let $V$ be a real vector bundle  on $M$ with structure group 
$Spin(2l)$. Similarly we can assume that $V$ has an $S^1$-equivariant spin 
structure without loss of generality.

The purpose of this part is to prove the elliptic operators introduced by 
Witten [{\bf W}] are rigid in the family case, at least at the equivariant 
Chern character level. Let us recall them more precisely.

For a vector bundle $E$ on $M$, let
\begin{eqnarray}\begin{array}{l}
S_t (E) = 1 + t E + t^2 S^2 E + \cdots ,\\
\Lambda_t (E) = 1 + tE + t^2 \Lambda^2 E + \cdots,
\end{array}\end{eqnarray}
be the symmetric and exterior power   operations in $K(M)[[t]]$.
Let 
\begin{eqnarray}\begin{array}{l}
\Theta'_q(TX) = \otimes_{n=1}^\infty \Lambda_{q^n} (TX) 
\otimes _{m=1}^\infty S_{q^m} (TX),\\
\Theta_q(TX) = \otimes_{n=1}^\infty \Lambda_{-q^{n-1/2}} (TX) 
\otimes _{m=1}^\infty S_{q^m} (TX),\\
\Theta_{-q}(TX) = \otimes_{n=1}^\infty \Lambda_{q^{n-1/2}} (TX) 
\otimes _{m=1}^\infty S_{q^m} (TX).
\end{array}\end{eqnarray}
We also define the following elements in $K(M) [[q^{1/2}]]$:
\begin{eqnarray}\begin{array}{l}
\Theta'_q(TX|V) = \otimes_{n=1}^\infty \Lambda_{q^n} (V) 
\otimes _{m=1}^\infty S_{q^m} (TX),\\
\Theta_q(TX|V) = \otimes_{n=1}^\infty \Lambda_{-q^{n-1/2}} (V) 
\otimes _{m=1}^\infty S_{q^m} (TX),\\
\Theta_{-q}(TX|V) = \otimes_{n=1}^\infty \Lambda_{q^{n-1/2}} (V) 
\otimes _{m=1}^\infty S_{q^m} (TX),\\
\Theta^*_q(TX|V) = \otimes_{n=1}^\infty \Lambda_{-q^n} (V) 
\otimes _{m=1}^\infty S_{q^m} (TX).
\end{array}\end{eqnarray}
Let $p_1(\cdot)_{S^1}$ 
denote the first $S^1$-equivariant Pontrjagin class and 
$\Delta(V) = \Delta^+ (V) \oplus \Delta ^- (V) $ be the spinor bundle of $V$.
 
In the following, we denote by $D^X\otimes W$ the Dirac operator 
on $\Delta \otimes W$ as defined in Section 1. We also write 
$d_s^X=D^X\otimes \Delta(TX)$. The following theorem is the family 
analogue of the Witten rigidity theorems as proved in [{\bf BT}], [{\bf T}]
 and [{\bf Liu2}].

\begin{thm}(a) The family operators $d_s^X\otimes \Theta'_q(TX)$, 
$D^X \otimes \Theta_q(TX)$ and $D^X \otimes \Theta_{-q}(TX) $ are rigid.

(b) If  $p_1(V)_{S^1} = p_1(TX)_{S^1}$, then 
$D^X \otimes\Delta(V) \otimes\Theta'_q(TX|V) $,
 $D^X \otimes (\Delta^+ (V) - \Delta ^- (V))\otimes \Theta^*_q(TX|V)$, 
 $D^X \otimes \Theta_q(TX|V) $ and $D^X \otimes \Theta_{-q}(TX|V)$ are rigid.
\end{thm}

\subsection{ \normalsize Proof of the family rigidity theorem}

For $\tau \in \bH = \{ \tau \in \bC; {\rm Im} \tau >0\}$,
 $q= e^{ 2\pi i \tau}$, let  
\begin{eqnarray}\begin{array}{l}
\theta_3(v, \tau)=c(q)\Pi_{n=1}^\infty (1 + q^{n-1/2} e^{2 \pi i v}) 
\Pi_{n=1}^\infty (1 + q^{n-1/2} e ^{-2 \pi i v}),\\
\theta_2(v, \tau)=c(q)\Pi_{n=1}^\infty (1 - q^{n-1/2} e^{2 \pi i v}) 
\Pi_{n=1}^\infty (1 - q^{n-1/2} e ^{-2 \pi i v}),\\
\theta_1(v, \tau)=c(q) q^{1/8} 2 \cos(\pi  v)
\Pi_{n=1}^\infty (1 + q^{n} e^{2 \pi i v}) 
\Pi_{n=1}^\infty (1 + q^{n} e ^{-2 \pi i v}),\\
\theta (v, \tau)=c(q)q^{1/8} 2 \sin (\pi v )
\Pi_{n=1}^\infty (1 - q^{n} e^{2 \pi i v}) 
\Pi_{n=1}^\infty (1 - q^{n} e ^{-2 \pi i v}).
\end{array}\end{eqnarray}
be the classical Jacobi theta functions [{\bf Ch}], 
where $c(q)= \Pi_{n=1}^\infty (1 - q^{n} )$.

Let $g= e^{2 \pi i t}\in S^1$ be a generator of the action group.
Let $\{M_{\alpha}\} $ be the fixed submanifolds of the circle action. 
Then $\pi: M_{\alpha}\to B$ be a submersion with fibre $X_{\alpha}$.
We have the following equivariant decomposition  of $TX$
\begin{eqnarray}
TX_{|M_{\alpha}} = N_1 \oplus \cdots \oplus N_h \oplus TX_{\alpha},
\end{eqnarray}
Here $N_{\gamma}$ is a complex vector bundle such that $g$ acts on it 
by $e^{2 \pi i m_{\gamma} t}$. We denote the Chern roots of
 $N_{\gamma}$ by $\{2 \pi i x_{\gamma} ^j\}$,  and the Chern roots of 
$TX_{\alpha} \otimes_{\sR} \bC$ by $\{ \pm 2 \pi i y_j\}$. 
We will write $\dim_{\sC} N_\gamma = d(m_\gamma)$
and $\dim X_\alpha = 2 k_\alpha$.

Similarly, let
\begin{eqnarray}
V_{|M_\alpha} = V_1 \oplus \cdots \oplus V_{l_0},
\end{eqnarray}
be the equivariant decomposition of $V$ restricted to $M_\alpha$. Assume that
$g$ acts on $V_v$ by $e^{2 \pi i n_v t}$, where some $n_v$ may be zero. 
We denote the Chern roots of $V_v$ by $2 \pi i u^j_v$. Let us write
$\dim_{\sR} V_v = 2d(n_v)$.

For $f(x)$ a holomorphic function, we denote by $f(y)(TX^g) = \Pi_j f(y_j)$, 
the symmetric polynomial which gives characteristic class of $TX^g$,
 and we use the same notation for $N_\gamma$. Now we define some functions on $\bC \times \bH$
 with values in $H^*(B)$,
\begin{eqnarray}\begin{array}{l}
\displaystyle{F_{d_s}(t, \tau) = \Sigma_{\alpha} \pi_* \Big [ 
\Big  (2 \pi  y  { \theta_1(y, \tau) \over \theta (y, \tau)}\Big )(TX^g) 
\Pi_{\gamma} \Big (i^{-1} { \theta_1(x_{\gamma} + m_{\gamma} t, \tau) 
\over \theta (x_\gamma + m_\gamma t, \tau)}\Big )(N_{\gamma}) \Big ],}\\
\displaystyle{F_D(t, \tau) = \Sigma_{\alpha} \pi_* \Big [ 
 \Big (2 \pi  y  { \theta_2(y, \tau) \over \theta (y, \tau)}\Big )(TX^g) 
\Pi_{\gamma} \Big (i^{-1} { \theta_2(x_{\gamma} + m_{\gamma} t, \tau) 
\over \theta (x_{\gamma} + m_{\gamma} t, \tau)}\Big )(N_{\gamma}) \Big ],}\\
\displaystyle{F_{-D}(t, \tau) = \Sigma_{\alpha} \pi_* \Big [ 
\Big  (2 \pi  y  { \theta_3(y, \tau) \over \theta (y, \tau)}\Big )(TX^g) 
\Pi_{\gamma} \Big (i^{-1} { \theta_3(x_{\gamma} + m_{\gamma} t, \tau) 
\over \theta (x_{\gamma} + m_{\gamma} t, \tau)}\Big )(N_{\gamma}) \Big ],}\\
\displaystyle{F^V_{d_s}(t, \tau) = i^{-k} \Sigma_{\alpha} \pi_* \Big [ 
\Big ({2 \pi i y \over \theta (y, \tau)}\Big )(TX^g)
 {\Pi_v  \theta_1(u_v + n_v t, \tau) (V_v) 
\over \Pi_{\gamma}\theta (x_{\gamma} + m_{\gamma} t, \tau)
(N_{\gamma})}\Big ],}\\
\displaystyle{F^V_{D}(t, \tau) = i^{-k} \Sigma_{\alpha} \pi_* \Big [ 
\Big ({2 \pi i y \over \theta (y, \tau)}\Big )(TX^g) 
 {\Pi_v  \theta_2(u_v + n_v t, \tau) (V_v) 
\over \Pi_{\gamma}\theta (x_{\gamma} + m_{\gamma} t, \tau)
(N_{\gamma})}\Big ],}\\
\displaystyle{F^V_{-D}(t, \tau) = i^{-k} \Sigma_{\alpha} \pi_* \Big [ 
\Big ({2 \pi i y \over \theta (y, \tau)}\Big )(TX^g) 
 {\Pi_v  \theta_3(u_v + n_v t, \tau) (V_v) 
\over \Pi_{\gamma}\theta (x_{\gamma} + m_{\gamma} t, \tau)
(N_{\gamma})}\Big ],}\\
\displaystyle{F^V_{D^*}(t, \tau) = i^{-k+ l} \Sigma_{\alpha} \pi_* \Big [ 
\Big ({2 \pi i y \over \theta (y, \tau)}\Big )(TX^g) 
 {\Pi_v  \theta(u_v + n_v t, \tau) (V_v) 
\over \Pi_{\gamma} \theta (x_{\gamma} + m_{\gamma} t, \tau)
(N_{\gamma})}\Big ].}
\end{array}\end{eqnarray}

By Theorem 1.1 and [{\bf LaM}, p238],  we get, 
for $t\in [0,1]\setminus \bQ$ and  $g=e^{2 \pi i t}$, 
\begin{eqnarray} \quad \begin{array}{l}
F_{d_s}(t, \tau)= \ch_g \Big ( {\rm Ind} (d_s^X \otimes \Theta'_q(TX))\Big ),\\
F_D(t, \tau) = q^{-k/8}\ch_g \Big ( 
{\rm Ind}(D^X \otimes \Theta_q(TX))\Big ),\\
F_{-D}(t, \tau) = q^{-k/8}\ch_g \Big ( {\rm Ind}
(D^X \otimes \Theta_{-q}(TX) )\Big ),\\
F^V_{d_s}(t, \tau)= c(q)^{l-k} q^{(l-k)/8}
\ch_g \Big ( {\rm Ind} (D^X \otimes\Delta(V) 
\otimes\Theta'_q(TX|V) )\Big ),\\
F^V_{D}(t, \tau)= c(q)^{l-k} q^{-k/8}\ch_g \Big ( {\rm Ind} 
(D^X \otimes \Theta_{q}(TX|V) )  \Big  )\\
F^V_{-D}(t, \tau) = c(q)^{l-k} q^{-k/8}\ch_g \Big ( {\rm Ind} (D^X \otimes 
\Theta_{-q}(TX|V))\Big ),\\
F^V_{D^*}(t, \tau)=(-1)^l c(q)^{l-k} q^{(l-k)/8}
\ch_g\Big  ({\rm Ind} (D^X \otimes 
(\Delta^+ (V) - \Delta ^- (V))\\
\hspace*{70mm} \otimes \Theta^*_q(TX|V))\Big ).
\end{array}\end{eqnarray}

Considered  as functions of $(t, \tau)$, we can obviously  extend these 
$F$'s and $F^V$'s to meromorphic functions on $\bC \times \bH$. Note that 
these functions are holomorphic in $\tau$. The rigidity theorems are equivalent to 
the statement that these  $F$'s and $F^V$'s
are independent of $t$. As explained in [{\bf Liu2}], 
we will prove it in two steps: i) we show that these $F,\ F^V$ are 
doubly periodic in $t$;  ii) we prove they are holomorphic in $t$. 
Then it is trivial to see that they are constant in $t$.

\begin{lemma} (a) For $a, b \in 2 \bZ$, $F_{d_s}(t, \tau), F_D(t, \tau)$ 
and $ F_{-D}(t, \tau)$ are invariant under the action
\begin{eqnarray}
U: t \to t + a \tau + b
\end{eqnarray}

(b) If $p_1(V)_{S^1} = p_1(TX)_{S^1}$, then $F^V_{d_s}(t, \tau), 
F^V_{D}(t, \tau)$, $F^V_{-D}(t, \tau)$ and  $F^V_{D^*}(t, \tau)$  
are invariant under $U$.
\end{lemma}

$Proof$: Recall  that we have the following transformation formulas of theta-functions [{\bf Ch}]:
\begin{eqnarray}\begin{array}{l}
\theta (t+1, \tau) = -\theta (t,\tau), \qquad 
\theta ( t+ \tau, \tau)= - q^{-1/2} e^{- 2 \pi i t} \theta (t, \tau),\\
\theta_1 (t+1, \tau) = -\theta_1 (t,\tau), \qquad 
\theta_1 ( t+ \tau, \tau)=  q^{-1/2} e^{- 2 \pi i t} \theta_1 (t, \tau),\\
\theta_2 (t+1, \tau) = \theta_2 (t,\tau), \qquad 
\theta_2 ( t+ \tau, \tau)= - q^{-1/2} e^{- 2 \pi i t} \theta_2 (t, \tau),\\
\theta_3 (t+1, \tau) = \theta_3 (t,\tau), \qquad 
\theta_3 ( t+ \tau, \tau)=  q^{-1/2} e^{- 2 \pi i t} \theta_3 (t, \tau).
\end{array}\end{eqnarray}
From these, for $\theta_v= \theta, \theta_1, \theta_2, \theta_3$ and  
$(a,b) \in (2 \bZ)^2$, $l\in \bZ$, we get
\begin{eqnarray}\begin{array}{l}
\theta_v (x + l(t+a\tau +b), \tau) =
e^{-\pi i (2 l a x + 2l^2 a t + l^2 a^2 \tau)} \theta_v (x + lt, \tau)
\end{array}\end{eqnarray}
which proves (a).

To prove (b), note that, since $p_1(V)_{S^1} = p_1(TX)_{S^1}$, we have
\begin{eqnarray}\begin{array}{l}
 \Sigma_{v,j} (u_v^j + n_v t)^2 
= \Sigma_j (y_j)^2 + \Sigma_{\gamma,j} (x_{\gamma}^j + m_{\gamma} t)^2.
\end{array}\end{eqnarray}
This implies the equalities:
\begin{eqnarray}\begin{array}{l}
\Sigma_{v,j} n_v u_v^j = \Sigma_{\gamma,j} m_{\gamma} x_{\gamma}^j,\\
\Sigma_{\gamma} m_{\gamma}^2 d(m_{\gamma}) = \Sigma_v n_v^2 d(n_v),
\end{array}\end{eqnarray}
which together with (2.11) proves (b).\hfill $\blacksquare$\\

For $g= \left ( \begin{array}{l} a \quad b\\
c \quad d \end{array} \right ) \in SL_2 (\bZ)$,  we define its modular 
transformation on $\bC \times \bH$ by
\begin{eqnarray}\begin{array}{l}
\displaystyle{
g(t, \tau) = \left ( { t \over c \tau + d}, {a \tau + b \over c \tau + d}
\right ). }
\end{array}\end{eqnarray}

The two generators of $SL_2(\bZ)$ are 
\begin{eqnarray} S = \left (  \begin{array}{l} 0 \quad -1 \\
1 \qquad    0
 \end{array} \right ),  \
T = \left (  \begin{array}{l} 1 \quad 1 \\
0 \quad 1
 \end{array} \right ),
\end{eqnarray}
which act on $\bC \times \bH$ in the following way:
\begin{eqnarray}
S(t,\tau) = \Big ({t \over \tau}, - {1 \over \tau}\Big ),\  
T(t, \tau) = (t, \tau +1).
\end{eqnarray}

Let $\Psi_{\tau}$ be the scaling homomorphism from $\Lambda (T^* B)$ into 
itself $: \beta \to \tau^{{1 \over 2} {\rm \scriptsize deg}\beta} \beta$.

\begin{lemma}(a) We have the following identities:
\begin{eqnarray}\begin{array}{l}
F_{d_s}({t \over \tau}, - {1 \over \tau}) 
= i^{k}\Psi_{\tau}F_D(t, \tau), \quad
F_{d_s}(t, \tau+1) = F_{d_s}(t, \tau),\\
 F_{-D}({t \over \tau}, - {1 \over \tau})
= i^{k} \Psi_{\tau}F_{-D}(t, \tau), \quad
F_{D}(t, \tau+1)= F_{-D}(t, \tau) e^{-{\pi i \over 4} k}.
\end{array}\end{eqnarray}

(b) If  $p_1(V)_{S^1} = p_1(TX)_{S^1}$, then we have 
\begin{eqnarray}\qquad \begin{array}{l}
F^V_{d_s}({t \over \tau}, - {1 \over \tau})= 
({\tau \over i})^{l- k \over 2} i^{k}\Psi_{\tau} F^V_{D}(t, \tau),\quad 
F^V_{d_s}(t, \tau + 1)= e^{-{\pi i \over 4}(k - l)} F^V_{d_s}(t, \tau),\\

 F^V_{-D}({t \over \tau}, - {1 \over \tau})= 
({\tau \over i})^{l- k \over 2} i^{k}\Psi_{\tau} F^V_{-D}(t, \tau),\quad 
 F^V_{D}(t, \tau + 1) = e^{- { \pi i \over 4} k} F^V_{-D}(t, \tau),\\

F^V_{D^*}({t \over \tau}, - {1 \over \tau}) = 
({\tau \over i})^{l- k \over 2} i^{k-l}\Psi_{\tau}
 F^V_{D^*}(t, \tau),\quad
F^V_{D^*}(t, \tau+1) = 
e^{-{\pi i \over 4}(k - l)} F^V_{D^*}(t, \tau).
\end{array}\end{eqnarray}
\end{lemma}

$Proof$: By [{\bf Ch}], we have the following transformation formulas 
for the Jacobi theta-functions:
\begin{eqnarray}\begin{array}{l}
\theta ({t \over \tau}, - {1 \over \tau})= {1 \over i} \sqrt{\tau \over i}
 e^{\pi i t^2 \over \tau} \theta (t,\tau), \quad 
 \theta (t, \tau+1) = e^{ \pi i \over 4} \theta (t, \tau),\\
\theta_1 ({t \over \tau}, - {1 \over \tau})=  \sqrt{\tau \over i}
 e^{\pi i t^2 \over \tau} \theta_2 (t,\tau), \quad 
 \theta_1 (t, \tau+1) = e^{ \pi i \over 4} \theta_1 (t, \tau),\\

\theta_2 ({t \over \tau}, - {1 \over \tau})=  \sqrt{\tau \over i}
 e^{\pi i t^2 \over \tau} \theta_1 (t,\tau), \quad 
 \theta_2 (t, \tau+1) =  \theta_3 (t, \tau),\\
\theta_3 ({t \over \tau}, - {1 \over \tau})=  \sqrt{\tau \over i}
 e^{\pi i t^2 \over \tau} \theta_3 (t,\tau), \quad 
 \theta_3 (t, \tau+1) =  \theta_2 (t, \tau).
\end{array}\end{eqnarray}

The action  of $T$ on the functions $F$ and $F^V$ are quite simple, and we 
leave the proof to the reader. Here we only check the action of $S$.
 By (2.19), we get 
\begin{eqnarray}\begin{array}{l}
\displaystyle{F_{d_s}({t \over \tau}, - {1 \over \tau})=
\Sigma_{\alpha} \pi_* \Big [ 
\Big  (2 \pi  y  { \theta_1(y,  - {1 \over \tau}) 
\over \theta (y,- {1 \over \tau} )}\Big  )(TX^g) 
\Pi_{\gamma} \Big (i^{-1} { \theta_1(x_{\gamma} 
+ m_{\gamma} {t \over \tau},  - {1 \over \tau}) 
\over \theta (x_{\gamma} + m_{\gamma} {t \over \tau},  
- {1 \over \tau}) }\Big )(N_{\gamma}) \Big ]  }\\
\displaystyle{\hspace*{10mm} 
= \Sigma_{\alpha}  i^{k} \tau ^{-k_\alpha}\pi_* \Big [ 
\Big (2 \pi  \tau y  { \theta_1(\tau y, \tau) 
\over \theta (\tau y, \tau)}\Big  )(TX^g) 
\Pi_{\gamma}\Big (i^{-1} { \theta_1(\tau x_{\gamma} + m_{\gamma} t, \tau) 
\over \theta (\tau x_{\gamma} + m_{\gamma}t, \tau)}\Big )(N_{\gamma}) \Big ] }
\end{array}\end{eqnarray}
If $ \alpha$ is a  differential form on $B$, we denote by $\{ \alpha\}^{(p)}$
the component of degree $p$ of $\alpha$. It is easy to see that (2.17)for $F_{d_s}$ follows from the following identity: 
\begin{eqnarray}\begin{array}{l}
\displaystyle{\tau ^{-k_\alpha}
\left \{\pi_* \Big [ 
\Big ( \tau y {\theta_1(\tau y, \tau) \over \theta (\tau y, \tau)}\Big )(TX^g) 
\Pi_{\gamma} \Big (i^{-1} { \theta_1(\tau x_{\gamma} + m_{\gamma} t, \tau) 
\over \theta (\tau x_{\gamma} + m_{\gamma} t, \tau)}\Big )(N_{\gamma})
 \Big ] \right \}^{(2p)} }\\
\displaystyle{
= \tau^{p} \left \{\pi_* \Big [\Big (  y  { \theta_1(y, \tau) 
\over \theta (y, \tau)}\Big )(TX^g) 
\Pi_{\gamma} \Big (i^{-1} { \theta_1(x_{\gamma} + m_{\gamma} t, \tau) 
\over \theta (x_{\gamma} + m_{\gamma} t, \tau)}\Big )(N_{\gamma}) \Big ]
 \right \}^{(2p)}}
\end{array}\end{eqnarray}
By looking at the degree $2(p+k_\alpha)$ part, that is the $(p+k_\alpha)$-th 
homogeneous terms of the polynomials in $x$'s and $y$'s,
 on both sides, we immediately get (2.21).

From (2.7), (2.20) and  (2.21), we obtain 
\begin{eqnarray}
\{F_{d_s}({t \over \tau}, - {1 \over \tau}) \}^{(2p)}
= i^{k} \tau^p \{F_D(t, \tau) \}^{(2p)},
\end{eqnarray}
which completes the proof of (2.17) for $F_{d_s}$. The other identities in (2.17) can be verified in the same way.

By using (2.12), (2.19) and the same trick as in the proof of (2.17), 
we can obtain the identities in (2.18). This completes the proof of Lemma 2.2. \hfill $\blacksquare$\\

The following lemma implies that the index theory comes in to cancel part of the 
poles of the functions $F$'s and $F^V$'s. 

\begin{lemma} If $TX$ and $V$ are spin, then all of the $F$'s and $F^V$'s 
above are holomorphic in $(t, \tau)$ for $(t, \tau) \in \bR \times \bH$.
\end{lemma}

$Proof$: Let $z= e^{2 \pi i t}$ and $l'=\dim M$. 
 We will consider these $F$'s and $F^V$'s as meromorphic functions of two
complex variables $(z, q)$ with values in $H^*(B)$. 

 i) Let $N= {\rm max}_{\alpha, \gamma} |m_\gamma|$. Denote by $D_N\subset \bC ^2$ the domain 
\begin{eqnarray}\begin{array}{l}
|q|^{1/N} < |z| < |q|^{-1/N}, 0< |q| < 1.
\end{array}\end{eqnarray}
Let $f_\alpha$ be the contribution of the component $M^\alpha$ in the functions
 $F$'s and $c(q)^{k-l}F^V$'s. Then in $D_N$, by (2.4), (2.7), 
it is easy to see that  $f_\alpha$ has expansions of the form
\begin{eqnarray}\begin{array}{l}
q^{-a/8} \Pi_\gamma (z^{m_\gamma} -1)^{-l' d(m_\gamma)} \Sigma_{n=1}^\infty
b_{\alpha, n}(z) q^n,
\end{array}\end{eqnarray}
where $a$ is an integer and  $h_\alpha (z,q)= \Sigma_{n=1}^\infty
b_{\alpha, n}(z) q^n$ is a  holomorphic function of $(z,q)\in D_N$, 
and $b_{\alpha, n}(z)$ are polynomial functions of $z$.
So as meromorphic functions, these  $F$'s and $c(q)^{k-l}F^V$'s 
have expansions of the form
\begin{eqnarray}\begin{array}{l}
q^{-a/8}  \Sigma_{n=1}^\infty
b_{ n}(z) q^n.
\end{array}\end{eqnarray}
with $b_n(z)$ rational function of $z$, which can only have poles 
on the unit circle $|z| = 1 \subset D_N$.

Now, if we multiply these  $F$'s and $c(q)^{k-l}F^V$'s by 
\begin{eqnarray}\begin{array}{l}
f(z)= \Pi_{\alpha} \Pi_\gamma (1 - z^{m_\gamma})^{l' d(m_\gamma)},
\end{array}\end{eqnarray}
we get holomorphic functions which have convergent power series expansions
 of the form
\begin{eqnarray}\begin{array}{l}
q^{-a/8} \Sigma_{n=1}^\infty c_n (z) q^n.
\end{array}\end{eqnarray}
with $\{ c_n (z)\}$  polynomial functions of $z$ in $D_N$. 

By comparing the above two expansions, we get for $n\in \bN$,
\begin{eqnarray}\begin{array}{l}
c_n(z) = f(z) b_n(z).
\end{array}\end{eqnarray}

ii) On the other hand, we can expand the Witten element $\Theta$'s into 
formal power series of the form $\Sigma_{n=0}^\infty R_n q^n$ with 
$R_n \in K(M)$. 
So, for $t\in [0, 1]\setminus 
 \bQ, z= e ^{2 \pi i t}$, we get a formal power serie of $q$ for 
these $F$'s and $c(q)^{k-l}F^V$'s:
\begin{eqnarray}\begin{array}{l}
q^{-a/8} \Sigma_{n=0}^{\infty} \ch_z ({\rm Ind}(D^X \otimes R_n)) q^n 
\end{array}\end{eqnarray}
with $a \in \bZ$.

By (1.6), we know
\begin{eqnarray}\begin{array}{l}
\ch_z ({\rm Ind}(D^X \otimes R_n))= \Sigma_{m= -N(n)}^{N(n)} a_{m,n} z^m.
\end{array}\end{eqnarray}
with $a_{m,n} \in H^*(B)$, 
and  $N(n)$ some positive integer  depending on $n$.

By comparing (2.7), (2.25) and (2.30), we get for $t\in [0, 1]\setminus 
 \bQ, z= e ^{2 \pi i t}$,
\begin{eqnarray}\begin{array}{l}
b_n (z) = \Sigma_{m= -N(n)}^{N(n)} a_{m,n} z^m.
\end{array}\end{eqnarray}
Since both sides are analytic functions of $z$, this equality holds 
for any $z\in \bC$.

By (2.28), (2.31), and the Weierstrass preparation theorem, we deduce that
\begin{eqnarray}
q^{-a/8}  \Sigma_{n=1}^\infty
b_{ n}(z) q^n = {1 \over f(z)} q^{-a/8} \Sigma_{n=1}^\infty
 c_n (z) q^n.
\end{eqnarray}
is holomorphic in $D_N$. 
Obviously $\bR \times \bH$ lies inside this domain. The proof of Lemma 2.3 is complete.\hfill $\blacksquare$\\

{\em Proof of the family rigidity theorem for spin manifolds}:
We will prove that these $F$'s and $F^V$'s are holomorphics on 
$\bC \times \bH$, which implies the rigidity theorem we want to prove.

We denote by $F$ one of the functions: $F$'s, $F^V$'s, $\Psi_{\tau}F$'s 
and $\Psi_{\tau}F^V$'s. From their expressions, we know the possible polar 
divisors of $F$ in $\bC \times \bH$ are of the form 
$t= {n \over l} (c\tau + d)$ 
with $n, c, d, l$ intergers and $(c,d)=1$ or $c=1$ and $d=0$.

We can always find intergers $a, b$ such that $ad-bc = 1$. Then the matrix 
$g=\left ( \begin{array}{l}d\quad -b\\
-c \quad a
\end{array} \right ) \in SL_2(\bZ)$ induces an action
\begin{eqnarray}\begin{array}{l}
\displaystyle{
F(g(t,\tau)) = F \Big ({t \over -c \tau + a},
 {d \tau -b \over -c \tau + a} \Big )}
\end{array}\end{eqnarray}
Now, if $t={n \over l} (c \tau + d)$ is a polar divisor of $F(t,\tau)$, 
then one polar divisor of $F(g(t,\tau))$ is given by
\begin{eqnarray}
{t \over -c \tau + a}= {n \over l}
 \Big ( c {d \tau -b \over -c \tau + a} + d \Big ),
\end{eqnarray}
which exactly gives $t = n/l$. 

But by Lemma 2.2, we know that, up to some constant, $F(g(t,\tau))$ is
still one of these $F$'s, $F^V$'s, $\Psi_{\tau}F$'s 
and $\Psi_{\tau}F^V$'s. This contradicts Lemma 2.3, therefore s completes the proof 
of Theorem 2.1.\hfill $\blacksquare$

\subsection{\normalsize A conjecture}

Motivated by the family rigidity theorem, Theorem 2.1, we and Zhang would like
 to make the following conjecture \\

{\bf Conjecture:} {\em The operators considered in Theorem 2.1 are rigid
 on the equivariant $K$-theory level.}\\

This means that, as elements in $K_G(B)$, the equivariant index bundles of those elliptic operators actually lie in $K(B)$. Note that this conjecture is more refined than Theorem 2.1, since the equivariant Chern character map is not an isomorphism. In [{\bf Z}], Zhang proved this for the canonical  $Spin^c$-Dirac operator on almost complex manifolds.

\section{ \normalsize Jacobi forms and vanishing theorems}
\setcounter{equation}{0}

In this Section, we generalize the rigidity theorems 
in the previous Section to the nonzero anomaly case, from which we derive a
 family of holomorphic Jacobi forms. As corollaries, we get many family 
vanishing theorems, especially a family $\widehat{ \frak U}$-vanishing 
theorem for loop space. This section  generalizes some results of 
[{\bf Liu4}, \S3] to the family case.

 This section is organized as follows: 
In Section 3.1, we state the generalization of  the rigidity theorems 
to the nonzero anomaly case. In Section 3.2, we prove this result. 
In Section 3.3, as corollaries, we derive several family vanishing theorems.

We will keep the notations of Section 2.

\subsection{ \normalsize  Nonzero anomaly}

Recall that the equivariant cohomology group $H^*_{S^1} (M, \bZ)$ 
of $M$ is defined by
\begin{eqnarray}
H^*_{S^1} (M, \bZ)= H^*(M \times_{S^1} ES^1, \bZ).
\end{eqnarray}
where $ES^1$ is the universal $S^1$-principal bundle over
the  classifying space $BS^1$ of $S^1$.
So $H^*_{S^1} (M, \bZ)$ is a module over $H^*(BS^1, \bZ)$ induced by the 
projection $\overline{\pi} : M\times _{S^1} ES^1\to BS^1$. 
Let $p_1(V)_{S^1}, p_1(TX)_{S^1} \in H^*_{S^1} (M, \bZ)$ be the equivariant
 first Pontrjagin classes of $V$ and $TX$ respectively. Also recall that 
\begin{eqnarray}
H^*(BS^1, \bZ)= \bZ [[u]]
\end{eqnarray}
with $u$ a generator of degree $2$.

In this section, we suppose that there exists some integer $n\in \bZ$ such that
\begin{eqnarray}
p_1(V)_{S^1}- p_1(TX)_{S^1} = n \cdot  \overline{\pi}^* u^2.
\end{eqnarray}
As in [{\bf Liu4}], we call $n$ the anomaly to rigidity.

Following [{\bf Liu2}], we introduce the following elements in $K(M)[[q^{1/2}]]$:
\begin{eqnarray}\begin{array}{l}
\Theta'_q(TX|V)_v = \otimes_{n=1}^\infty \Lambda_{q^n} (V- \dim V) 
\otimes _{m=1}^\infty S_{q^m} (TX- \dim X),\\
\Theta_q(TX|V)_v = \otimes_{n=1}^\infty \Lambda_{-q^{n-1/2}} (V- \dim V) 
\otimes _{m=1}^\infty S_{q^m} (TX- \dim X),\\
\Theta_{-q}(TX|V)_v = \otimes_{n=1}^\infty \Lambda_{q^{n-1/2}} (V- \dim V) 
\otimes _{m=1}^\infty S_{q^m} (TX- \dim X),\\
\Theta^*_q(TX|V)_v = \otimes_{n=1}^\infty \Lambda_{-q^n} (V- \dim V) 
\otimes _{m=1}^\infty S_{q^m} (TX- \dim X).
\end{array}\end{eqnarray}

For $g=e^{2 \pi i t},  q= e ^{2 \pi i \tau}$, 
with $ (t, \tau)\in \bR \times \bH$, we denote the equivariant Chern character 
of the index bundle  of $D^X \otimes\Delta(V) \otimes\Theta'_q(TX|V)_v $,
 $D^X \otimes \Theta_q(TX|V)_v$, $D^X \otimes \Theta_{-q}(TX|V)_v$, 
and $D^X \otimes (\Delta^+ (V) - \Delta ^- (V))\otimes \Theta^*_q(TX|V)_v$
by $2^l F^V_{d_s,v}(t, \tau),F^V_{D,v}(t, \tau)$, $F^V_{-D,v}(t, \tau)$ and 
 $(-1)^l F^V_{D^*,v}(t, \tau)$ respectively.
Similarly, we denote by $H(t, \tau)$  the  equivariant Chern character 
of the index bundle  of 
\begin{eqnarray}
D\otimes \otimes _{m=1}^\infty S_{q^m} (TX- \dim X). \nonumber 
\end{eqnarray}
Later we will consider these functions as the extensions of these  functions
from the unit  circle with variable $e^{2\pi it}$ to the complex plane with values in $H^*(B)$. For $\alpha$ a differential form on $B$, we denote by $\{\alpha\}^{(p)}$ the degree $p$ component of $\alpha$.

The purpose of this section is to prove the following theorem 
which generalizes the family rigidity  theorems to the nonzero anomaly case.

\begin{thm} Assume $p_1(V)_{S^1}- p_1(TX)_{S^1} = n\cdot \overline{\pi}^* u^2$ 
with $n\in \bZ$. Then for $p\in \bN$, 
$\{F^V_{d_s,v}(t, \tau)\}^{(2p)}$, $\{F^V_{D,v}(t, \tau)\}^{(2p)}$,
 $ \{F^V_{-D,v}(t, \tau)\}^{(2p)}$
 are holomorpic Jacobi forms of index $n/2$ and weight $k+p$ over
 $(2\bZ)^2 \rtimes \Gamma$ with $\Gamma$ equal to $\Gamma_0(2), \Gamma^0(2), 
\Gamma_\theta$ respectively, and $\{F^V_{D^*,v}(t, \tau)\}^{(2p)}$ 
is a holomorphic  Jacobi form of index $n/2$ 
and weight $k-l+p$ over $(2 \bZ)^2 \rtimes SL_2(\bZ)$.
\end{thm}

See Section 3.2 for the definitions of these modular subgroups, 
$\Gamma_0(2), \Gamma^0(2)$ and  $\Gamma_\theta$. 

\subsection{\normalsize Proof of Theorem 3.1}

Recall that a (meromorphic) Jacobi form  of index $m$ and weight $l$ over 
$L\rtimes \Gamma$, where $L$ is an integral lattice  in the complex plane 
$\bC$ preserved by the modular subgroup $\Gamma \subset SL_2(\bZ)$,
 is a (meromorphic) function $F(t, \tau)$ on $\bC \times \bH$ such that
\begin{eqnarray}\begin{array}{l}
\displaystyle{F({t \over c \tau + d}, { a \tau + b \over c \tau +d}) 
= (c \tau + d)^l e^{2 \pi i m ( c t^2 / (c \tau + d))} F(t, \tau),}\\
\displaystyle{F(t + \lambda \tau + \mu, \tau) = 
e ^{- 2 \pi i m ( \lambda ^2 \tau + 2 \lambda t)} F(t, \tau),}
\end{array}\end{eqnarray}
where $(\lambda, \mu) \in L$, and $g = \left (\begin{array}{l} a \quad b\\
c \quad d
\end{array}  \right ) \in \Gamma$. If $F$ is holomorphic on $\bC \times \bH$, 
we say that $F$ is a holomorphic Jacobi form.

Now, we start to prove Theorem 3.1.  Let $g= e^{2 \pi i t} \in S^1$ 
be a generator of  the action group. For $\alpha = 1,2,3$, let
\begin{eqnarray}
\theta'(0, \tau) = {\partial \over \partial t} \theta(t, \tau)|_{t=0},
\qquad \theta_{\alpha}(0, \tau) = \theta_{\alpha}(t, \tau)|_{t=0}
\end{eqnarray}

By applying Theorem 1.1, we get
\begin{eqnarray}\begin{array}{l}
\displaystyle{F^V_{d_s,v}(t, \tau) = (2 \pi)^{-k}
{ \theta'(0, \tau)^k \over \theta_1(0, \tau)^l}
F^V_{d_s}(t, \tau), } \\
\displaystyle{F^V_{D,v}(t, \tau) = (2 \pi)^{-k}
{ \theta'(0, \tau)^k \over \theta_2(0, \tau)^l}F^V_{D}(t, \tau), }\\

\displaystyle{F^V_{-D,v}(t, \tau) = (2 \pi)^{-k}
{ \theta'(0, \tau)^k \over \theta_3(0, \tau)^l}
F^V_{-D}(t, \tau) ,}\\
\displaystyle{F^V_{D^*,v}(t, \tau) =  (2 \pi)^{l-k}
 \theta'(0, \tau)^{k-l}
F^V_{D^*}(t, \tau),}\\
\displaystyle{H(t,\tau) = (2 \pi i)^{-k}\Sigma_{\alpha} \pi_* \Big [ 
\Big ({2 \pi i y \over \theta (y, \tau)}\Big )(TX^g) 
 { \theta'(0, \tau)^k
\over \Pi_{\gamma}\theta (x_{\gamma} + m_{\gamma} t, \tau)
(N_{\gamma})}\Big ].}
\end{array}\end{eqnarray}

\begin{lemma} If $p_1(V)_{S^1}- p_1(TX)_{S^1} = n \cdot \overline{\pi}^* u^2$,
 we have 
\begin{eqnarray}\qquad \begin{array}{l}
F^V_{d_s,v}({t \over \tau}, - {1 \over \tau})= 
\tau^k e^{\pi i n t^2/\tau} \Psi_{\tau} F^V_{D,v}(t, \tau),\quad 
F^V_{d_s,v}(t, \tau + 1)=  F^V_{d_s,v}(t, \tau),\\

 F^V_{-D,v}({t \over \tau}, - {1 \over \tau})= 
\tau^k e^{\pi i n t^2/\tau} \Psi_{\tau} F^V_{-D,v}(t, \tau),\quad 
 F^V_{D,v}(t, \tau + 1) =  F^V_{-D,v}(t, \tau),\\

F^V_{D^*,v}({t \over \tau}, - {1 \over \tau}) = 
\tau^{k-l} e^{\pi i n t^2/\tau} \Psi_{\tau}
 F^V_{D^*,v}(t, \tau),\quad
F^V_{D^*,v}(t, \tau+1) =   F^V_{D^*,v}(t, \tau).
\end{array}\end{eqnarray}
If $p_1(TX)_{S^1} = -n \cdot \overline{\pi}^* u^2$, then
\begin{eqnarray}
H({t \over \tau}, - {1 \over \tau})= 
\tau^k e^{\pi i n t^2/\tau} \Psi_{\tau} H(t, \tau), \quad
H(t,\tau+1) =   H(t, \tau).  
\end{eqnarray}
\end{lemma}

$Proof$: First recall that the condition on the first equivariant 
Pontrjagin classes implies the equality
\begin{eqnarray}
\Sigma_{v,j} (u_v^j + n_v t)^2 -\Big ( \Sigma_j (y_j)^2 
+ \Sigma_{\gamma,j} (x_{\gamma}^j + m_{\gamma} t)^2\Big ) = n\cdot t^2
\end{eqnarray}
which gives the equalities:
\begin{eqnarray}\begin{array}{l}
\Sigma_v n_v^2 d(n_v)-\Sigma_{\gamma} m_{\gamma}^2 d(m_{\gamma}) =n, \quad 
\Sigma_{v,j} n_v u_v^j = \Sigma_{\gamma,j} m_{\gamma} x_{\gamma}^j,\\
\Sigma_{v,j} (u_v^j)^2 =
\Sigma_j (y_j)^2 + \Sigma_{\gamma,j} (x_{\gamma}^j)^2.
\end{array}\end{eqnarray}
The action of $T$ on the functions $F$ and $F^V$ are quite easy to check and
 we leave the detail to the reader. We only check the action of $S$. By (2.7), 
(2.19), (3.7) and  (3.11), we have
\begin{eqnarray}  \qquad\begin{array}{l}
\displaystyle{F^V_{d_s,v}({t \over \tau}, - {1 \over \tau})= (2 \pi i)^{-k} }\\
\displaystyle{\hspace*{20mm} 
 \Sigma_{\alpha} \pi_* \Big [ 
{ \theta'(0,- {1 \over \tau} )^k \over \theta_1(0, - {1 \over \tau})^l}
\Big ({2 \pi i y \over \theta (y, - {1 \over \tau})}\Big )(TX^g)
 {\Pi_v  \theta_1(u_v + n_v {t \over \tau}, - {1 \over \tau}) (V_v) 
\over \Pi_{\gamma}\theta (x_{\gamma} + m_{\gamma} {t \over \tau},
 - {1 \over \tau})(N_{\gamma})}\Big ]    }\\
\displaystyle{
= (2 \pi i)^{-k}\tau^k e^{\pi i (n t^2 /\tau)} \Sigma_{\alpha} \pi_* \Big [ 
{ \theta'(0, \tau)^k \over \theta_1(0, \tau)^l}
\Big ({2 \pi i y \over \theta (\tau y, \tau)}\Big )(TX^g)
 {\Pi_v  \theta_1(\tau u_v + n_v t, \tau) (V_v) 
\over \Pi_{\gamma}\theta (\tau x_{\gamma} + m_{\gamma} t, \tau)
(N_{\gamma})}\Big ]}.
 \end{array}\end{eqnarray}
As in (2.21), by comparing the $(p+k_\alpha)$-th homogeneous terms of 
the polynomials in $x$'s,  $y$'s and $u$'s on both side, 
we find the following equation
\begin{eqnarray}\begin{array}{l}
\displaystyle{ \Big \{ \pi_* \Big [ 
\Big ({2 \pi i y \over \theta (\tau y, \tau)}\Big )(TX^g)
 {\Pi_v  \theta_1(\tau u_v + n_v t, \tau) (V_v) 
\over \Pi_{\gamma}\theta (\tau x_{\gamma} + m_{\gamma} t, \tau)
(N_{\gamma})}\Big ] \Big \}^{(2p)} }\\
\displaystyle{ = \tau^p  \pi_* \Big [ 
\Big ({2 \pi i y \over \theta ( y, \tau)}\Big )(TX^g)
 {\Pi_v  \theta_1( u_v + n_v t, \tau) (V_v) 
\over \Pi_{\gamma}\theta (x_{\gamma} + m_{\gamma} t, \tau)
(N_{\gamma})}\Big ] \Big \}^{(2p)} . }
\end{array}\end{eqnarray}
By (3.12), (3.13), we get the equation (3.8) for $F^V_{d_s, v}$.
We leave the other cases to the reader. \hfill $\blacksquare$\\

Recall the three modular subgroups:
\begin{eqnarray} \qquad \begin{array}{l}
\Gamma_0(2) = \left \{ \left ( \begin{array}{l} a \quad b\\
c \quad d  \end{array}\right ) \in SL_2 (\bZ)| c \equiv 0 (\rm mod 2)
 \right \}, \\
\Gamma^0(2) = \left \{ \left ( \begin{array}{l} a \quad b\\
c \quad d \end{array}\right ) \in SL_2 (\bZ)| b \equiv 0 (\rm mod 2) 
\right \},\\
\Gamma_{\theta} = \left \{ \left ( \begin{array}{l} a \quad b\\
c \quad d \end{array} \right ) \in SL_2 (\bZ)| 
\left ( \begin{array}{l} a \quad b\\
c \quad d \end{array} \right )
 \equiv   \left ( \begin{array}{l} 1 \quad 0\\
0 \quad 1 \end{array} \right ) {\rm  or } 
\left ( \begin{array}{l} 0 \quad 1 \\
1 \quad 0 \end{array} \right )(\rm mod 2) 
\right \}.
\end{array}\end{eqnarray}

\begin{lemma} If $p_1(V)_{S^1}- p_1(TX)_{S^1} = n \cdot \overline{\pi}^* u^2$, then
for $p\in \bN$, $\{F^V_{d_s,v}(t, \tau)\}^{(2p)}$ 
is a Jacobi form over $(2\bZ)^2 \rtimes \Gamma_0(2)$; $\{F^V_{D,v}(t, \tau)\}^{(2p)}$ is a Jacobi form over $(2\bZ)^2 \rtimes \Gamma^0(2)$; $\{F^V_{-D,v}(t, \tau)\}^{(2p)}$ is a Jacobi form over $(2\bZ)^2 \rtimes \Gamma_\theta$. If $p_1(TX)_{S^1}= -n \overline{\pi}^* u^2$, then $\{H(t, \tau)\}^{(2p)}$ 
is a Jacobi form over $(2\bZ)^2 \rtimes SL_2(\bZ)$. All of them are of index ${n \over 2}$ and weight $k+p$. 

The function $\{F^V_{D^*,v}(t, \tau)\}^{(2p)}$ is  a Jacobi form of index ${n \over 2}$ 
and weight $k-l+p$ over $(2\bZ)^2 \rtimes SL_2(\bZ)$.
\end{lemma}

$Proof$: By (2.19), (3.7), we know that these $F^V$'s and $H$ 
satisfy the second equation of the definition of Jacobi forms (3.5).

Recall that $T$ and $ST^2 ST$ generate $\Gamma_0(2)$, and  also $\Gamma^0(2)$ 
and $\Gamma_{\theta}$ are conjugate to $\Gamma_0(2)$ by $S$ and $TS$ 
respectively. By Lemma 3.1, and the above discussion, for $F^V$'s and $H$,
  we easily get the first equation  of (3.5). \hfill $\blacksquare$\\

For $g=  \left ( \begin{array}{l} a \quad b\\
c \quad d \end{array}\right ) \in SL_2 (\bZ)$, let us use the notation
\begin{eqnarray}
F(g(t,\tau))|_{m,l} = (c \tau + d)^{-l} e^{- 2 \pi i m c t^2/(c \tau +d)} 
F({t \over c \tau + d}, { a \tau +b \over  c\tau +d }).
\end{eqnarray}
to denote the action of $g$ on a Jacobi form $F$ of index $m$ and weight $l$.

By Lemma 3.1, for any function in $ F\in \{ {\{F^V\}^{(2p)}}, H^{(2p)} \}$, 
its modular  transformation $\{F\}^{(2p)}(g(t, \tau))|_{{n \over 2}, k+p}$ 
 (or $\{F\}^{(2p)}(g(t, \tau))|_{{n \over 2}, k-l+p}$)
is still one of the $\{ {\{F^V\}^{(2p)}}\}$'s and $H^{(2p)}$.
 Similar to   Lemma 2.3,  we have 

\begin{lemma} For any $g\in SL_2(\bZ)$, let $F(t, \tau)$ be one of the 
$\{F^V\}^{(2p)}$'s or $H^{(2p)}$, Then $F(g(t, \tau))|_{{n \over 2}, k+p}$
is holomorphic in $(t, \tau)$ for $t\in \bR$ and $\tau \in \bH$.
\end{lemma}
As in Lemma 2.3, this is the place where the index theory comes in to cancel part 
of the poles of these  functions. Of course, to use the index theory, 
we must use the spin conditions on $TX$ and $V$.

Now, we recall the following result  [{\bf Liu4}, Lemma 3.4]:
\begin{lemma} For a (meromorphic) Jacobi form $F(t, \tau)$ of index $m$ and 
weight $k$ over $L\rtimes \Gamma$, assume that $F$ may only have polar 
divisors of the form $t=(c \tau + d)/l$ in $\bC \times \bH$ for some 
integers $c,d$ and $l\neq 0$. If $F(g(t, \tau))|_{m,k}$ is holomorphic 
for $t\in \bR$, $\tau \in \bH$ for every $g\in SL_2(\bZ)$, 
then $F(t, \tau)$ is holomorphic for any $t\in \bC$ and $\tau \in \bH$.
\end{lemma}

{\em Proof of Theorem 3.1}: By Lemmas 3.1, 3.2, 3.3, we know that the 
$\{F^V\}^{(2p)}$'s and $H^{(2p)}$ satisfy the assumptions of Lemma 3.4. 
In fact, all of their  possible polar divisors are of the form 
$l= (c \tau + d)/m$ where $c,d$ are integers and $m$ is one of 
the exponents $\{ m_j\}$. The proof of Theorem 3.1 is complete.
\hfill $\blacksquare$

\subsection{\normalsize Family vanishing theorems for loop space}

The following lemma is established in [{\bf EZ}, Theorem 1.2]:

\begin{lemma} Let $F$ be a holomorphic Jacobi form of index $m$ and weight $k$.
 Then for fixed $\tau$, $F(t,\tau)$, if not identically zero, has exactly 
 $2m$ zeros in any fundamental domain for the action of the lattice on $\bC$.
\end{lemma}
This tells us that there are no holomorphic Jacobi forms of negative index.
 Therefore, if $m<0$, $F$ must be identically zero. 
If $m=0$, it is easy to see  that $F$ must be independent of $t$.

Combining Lemma 3.5 with Theorem 3.1, we have the following result.
\begin{cor} Let $M,B,V$ and $n$ be as in Theorem 3.1, If $n=0$, the 
equivariant Chern characters of the index bundle of  
 the elliptic operators in Theorem 3.1 are independent of $g\in S^1$. 
If $n<0$, then these equivariant Chern characters  are identically zero; 
in particular, 
the Chern character of the index bundle of these elliptic operators are zero.
\end{cor}

Another quite interesting consequence of the above discussions is 
the following family 
$\widehat{\frak U}$-vanishing theorem for loop space.

\begin{thm} Assume $M$ is connected and the $S^1$-action is nontrivial. 
If $p_1(TX)_{S^1} = n \cdot \overline{\pi}^* u^2$ for some 
integer $n$, then the equivariant Chern character of the index bundle,
 especially the Chern character of the index bundle, of 
$D\otimes \otimes _{m=1}^\infty S_{q^m} (TX- \dim X)$ is identically zero. 
\end{thm}

$Proof$: In fact, by (3.11), we know that 
\begin{eqnarray}
\Sigma_j m^2_j d(m_j)=n.
 \end{eqnarray}
So  the case $n<0$ can never happen. If $n=0$, then all the exponents 
$\{m_j\}$ are zero, so the $S^1$-action cannot  have a fixed point. 
By (2.7)and (3.7), we know that $H(t, \tau)$ is zero. 
For $n>0$, one  can apply Lemmas 
3.1, 3.4 and 3.5 to get the result.\hfill $\blacksquare$\\

As remarked in [{\bf Liu4}],  the fact that the index of 
$D\otimes \otimes _{m=1}^\infty S_{q^m} (TX- \dim X)$ 
is zero may be viewed as a loop space analogue of the famous 
$\widehat{\frak U}$-vanishing theorem of Atiyah and Hirzebruch [{\bf AH}]
for compact connected spin manifolds with non-trivial $S^1$-action. The reason is that, this operator corresponds to the Dirac operator on loop space $LX$, while the condition on $p_1(TX)_{S^1}$ is a condition for the existence of an equivariant spin structure on  $LX$. This property is one of  the most interesting and surprising properties of loop space.
Now, under the condition of Theorem 3.2, a very interesting question is to know
 when the index bundle of this elliptic operator is zero in $K(B)$.\\

{\bf Acknowledgements.}  We would like to thank W. Zhang for many interesting 
discussions. Part of this work was done while the second author was visiting 
IHES. He would like to thank Professor J. P. Bourguignon and IHES
for their hospitality.

\centerline{------------------------}
\vskip 6mm

Kefeng LIU,
Department of Mathematics, Stanford University, Stanford, CA 94305, USA.

{\em E-mail address}: kefeng@math.stanford.edu

\vskip 6mm

Xiaonan MA,
Humboldt-Universitat zu Berlin, Institut f\"ur Mathematik,unter den Linden 6,
D-10099 Berlin, Germany.

{\em E-mail address}: xiaonan@mathematik.hu-berlin.de

\end{document}